\documentclass{amsproc}

\begin{document}

\title{Branch cuts: writing, editing, and ramified complexities}
\author{Ursula Whitcher}
\address{Mathematical Reviews, American Mathematical Society, Ann Arbor, MI 48103}
\curraddr{}
\email{uaw@umich.edu}
\thanks{Thank you to the anonymous referees who pushed me to be bolder and more definite, the friends and family members who helped me think through where and why ambiguity was important to me, and all the people who have enthused with me about mathematics over the years.}
\date{}

\begin{abstract}
As I was preparing my tenure application, the University of Wisconsin Board of Regents voted to redefine tenure, removing many of the institution's historical protections. Reevaluating my career priorities in light of these changes and a resurgent two-body problem, I recognized that my fundamental goal was communicating mathematical ideas. I found a new role as an editor at Mathematical Reviews, part of the American Mathematical Society. To my surprise, thinking more about my identity as a writer and editor also changed my perspective on my own sexuality and gender identity, inspiring new approaches to leadership.
\end{abstract}

\maketitle

%* Going up for tenure as the Board of Regents voted to redefine tenure

\section{The Wisconsin Idea}

In the spring of 2011, between one job interview and another, I visited friends in Madison, Wisconsin. Downtown was thick with protesters, bundled up in scarves and hats. They were protesting the governor Scott Walker's decision to erase bargaining rights for public employee unions, a change that hit public schools and universities hard. My friends told me that a local pizza chain was delivering to protesters, and it had received orders from as far away as Syria: protesters in the midst of the Arab Spring were reaching out to support protesters here.

I knew, therefore, when I accepted a tenure-track job at the University of Wisconsin--Eau Claire, that Wisconsin universities were under attack. But I also knew that they had many defenders. Eau Claire felt, for reasons both concrete and ineffable, like a place where I could be useful, and could thrive. Professionally, I appreciated the fact that UWEC had a long-standing commitment to undergraduate research, including funds to pay student researchers for their work during the academic year. On a campus where many students supported themselves with a patchwork of loans and part-time jobs, the stipend was key to recruiting and engaging student researchers. I also had personal reasons to appreciate the location: Eau Claire was next door to Chippewa Falls, the home of Cray Supercomputing, and my spouse had interned for Cray as an engineering student in Seattle. As we were learning, serious two-body problems where only one partner is an academic are surprisingly rare. It's easier to find a job as an electrical engineer than as a pure mathematician, but the geographic constraints are different, and university colleagues rarely know how to help. The existing industry near Eau Claire seemed like rare good fortune.

My first semester was rough. My spouse was still working in California, where I had done a postdoc, and every time I flew back to visit him, I got sick. After Thanksgiving, I missed a connecting flight and spent a night in the Phoenix airport, too exhausted to negotiate the problem of leaving security and finding a hotel. 

I was teaching college algebra for the first time. I'd borrowed materials and class format from a loud and forthright male colleague, and his jovial, pushy classroom style didn't translate to my own more formal and self-effacing classroom persona. I was adjusting expectations from a tiny engineering college where students prided themselves on academic overwork to an institution where classes were only one obligation among many. 

But the hardest thing for me about teaching college algebra was that it reminded me of middle school. My family moved to the Portland, Oregon suburbs just before I started sixth grade, and I found myself in a newly built school in the middle of empty fields at the leading edge of suburbia. My new classmates informed me, in ways both obvious and subtle, that I was bad at being a girl. I was wearing the wrong clothes, or the right brands in the wrong ways. I hadn't realized I should wear a bra and shave my legs, I didn't know anything about the music on the radio, and most obvious of all, my hair was wrong---to start with, I wasn't blonde. I only escaped overtly homophobic bullying out of sheer na\"{i}vet\'e. I remember the most beautiful girl in my seventh-grade English class asking whether I had ever \emph{experienced} the feeling of attraction, and answering doggedly, ``No, never, not at all," while all the time noticing the way her hair fell against her cheek. Taking advanced math classes in middle school, and in particular the math team elective, wrecked the rest of my schedule and gave me an escape: I transferred out of that English class, I spent entire terms away from P.E. and the empty, echoing locker room, and I got to hang around eighth graders, who were too tall and grand to care what I was wearing. 

Here I was, then, at the beginning of another move to a new city, in a room full of blond Wisconsin students, the only woman without makeup, trying to teach the algebra I'd learned in seventh grade. I did my best to squash the associations. But when you're working really hard not to feel a feeling, the place that you end up is strange and flat, and other people make their own guesses about where the strangeness comes from. My students had their own baggage about algebra class, and no reason to be kind, or to anticipate kindness.

The bright spot in that first semester was my research students, who were rapidly figuring out how to build new four-dimensional polytopes. The following semester, they won the provost's research prize. That spring was good in other ways, too. I taught a cheerful business calculus course. To my astonishment, two of those students recruited me as the statutory advisor to the Equestrian Club. I signed their forms, in exchange for photographs of horses. My spouse and I bought a rambling Victorian house with a Persian rug set into the living-room floor. Wisconsin springs are snowy, but I got my own pair of snowshoes and joined my colleagues in exploring frozen streams and snow-covered pines.

Things really clicked into place when I started teaching Math for Liberal Arts. This was a 100-level class that wasn't a prerequisite for anything else, so I had complete freedom: my only charge was to persuade my students that math was bigger and more fascinating than they had realized. We talked about cryptography, voting theory, and the fourth dimension. We used measurements of gerrymandering to explore topological concepts such as connectedness; my students made brochures advertising unsolved math problems. In this context, being a little bit weird felt natural: I could be warm and enthusiastic about strange things in a genuine way.

Meanwhile, my research was going swimmingly. I wrote papers with my students bridging geometry and combinatorics. I had also formed a new collaboration focused on the intersections between mathematical physics, algebraic geometry, and number theory. I wrote grants to support the collaboration, and we met in Banff and San Jos\'e, scribbling correspondences and keeping an eye out for passing elk.

The University of Wisconsin System allows professors to split promotion and tenure, so I successfully applied for promotion to Associate Professor in my fourth year at UWEC, using credit from my time as a postdoc. 

\section{All paths to victory}

%* Tension between teaching and research--is it real? What are the gendered expectations? How did this affect my "ideal" job?
%* Applying to the Math Reviews job...

In my fifth year at UWEC, I was preparing my tenure portfolio, Scott Walker had survived a recall and been re-elected, and Wisconsin universities were again under attack. The target this time was the definition of tenure, which had previously been enshrined in the Wisconsin state constitution. By redefining tenure, the governor's hand-picked Board of Regents made it easier to remove entire departments by citing manufactured budgetary problems. 

I wasn't worried that the math department would disappear (though subsequently some University of Wisconsin campuses did consider worrisome program changes, potentially dropping majors such as history and Spanish). But I was worried about what the attack on tenure said about the state's attitude toward its university system, and what that would mean for my day-to-day job experience. Class sizes were getting bigger, and offering low-enrollment upper-division classes was getting harder. I loved the freedom that being a professor gave me to learn new things, but as administrative duties and the sheer volume of grading accumulated, that freedom seemed more and more hypothetical.

Meanwhile, the solution to my two-body problem was unstable. My spouse had indeed found a job in Chippewa Falls, with a supercomputing company that had spun off from Cray. But their supervisor had survived many years of corporate reconfigurations by being ruthlessly territorial: he could not admit that anyone else had expertise. The misery and frustration accumulated, in a kind of haze. 

Thus, as I submitted my tenure portfolio, I quietly applied for new jobs. As I did so, I thought carefully about the balance of teaching and research. All math-professor jobs involve some of each, but positions where teaching and research are equally balanced are rare. Most professorial positions expect, explicitly or tacitly, that the successful candidate will focus time and attention on one side or the other. 

Moreover, expectations about teaching and research are bound up with expectations about gender in complicated ways. In one paradigmatic story, a brilliant young man excels at mathematics, chess, and playing the piano, wins a number of math competitions, and becomes an acclaimed researcher. In another paradigmatic story, a girl who has always loved school and dreamed of being a teacher gets to college, realizes she could teach university students, and sets her sights on a PhD. As a professor, she nurtures many students, helping them realize their own potential. My own teaching persona is less overtly nurturing and more nerd bubbling over with obscure historical context. I knew from past job searches that I would mesh with some departmental cultures, and slam head-on into unspoken expectations in others. For example, I am not the right person to organize outreach programs for middle school girls!

Meanwhile, I found one strange job ad that didn't involve teaching at all. Mathematical Reviews (perhaps better known for its database MathSciNet), part of the American Mathematical Society, was advertising for a mathematician with expertise in algebraic geometry or group theory who had broad mathematical interests, strong writing skills, and significant experience beyond the PhD. I was used to broad calls for motivated teachers and researchers of every stripe. Here, in contrast, was a job that seemed targeted specifically at me. My research in algebraic geometry focused on a web of connections to number theory and theoretical physics; I was interested in expository writing and had won the MAA's Merten Hasse Prize for exposition; I was at the right career stage.

I interviewed in Ann Arbor in the dregs of a March snowstorm. The interview was subdued but friendly. I looked over my predecessor's shoulder as he assigned a paper for review that was directly relevant to my own current research project; I hand-edited a sample research description; I chatted with copyeditors and learned about the building's history as a pre-Prohibition brewery. I came away with freshly sharpened green and purple editing pencils, a strong desire to spend more time in Ann Arbor coffee shops, and quiet curiosity. I flew straight into a snowstorm in O'Hare.

I got a phone call offering me the Math Reviews job on spring break, while I was at a conference in North Carolina. The sky was brilliant blue, the tulip trees were blooming, and the salary was a huge percentage increase over my state university salary. I called the most competitive of my several beloved postdoctoral mentors (the one who had gone to MIT and used to coach the Putnam) and got him to walk me through negotiating anyway. He pointed out that I was comparing a nine-month salary to a twelve-month salary; it was still a significant percentage increase. I screwed up my courage, for feminism, and talked my way into a signing bonus.

I still had to answer the bigger question: was I really going to walk away from a job as a tenured professor? There's a huge mystique bound up in tenure, but fundamentally it's a job recruitment tool: it offers security and status. I had already earned the status. The security associated with tenure in the University of Wisconsin System felt more and more like a fiction. 

Thinking about leaving a teaching role was stranger. I asked myself what I most valued about teaching mathematics, and realized it was being part of a community that valued mathematical communication. I wanted to talk about mathematics, and to share it. Formal teaching was one way to do that, but definitely not the only way.

My spouse and I found a new solution to our two-body problem before we officially arrived in Michigan, through a fortuitous combination of training and circumstance. We had a high school friend who was just starting a tenure-track position as an experimental physicist, and her department needed an electrical engineer. Thus, my spouse started a new role as an independent contractor, splitting their time between design for the IceCube neutrino experiment and the Large Hadron Collider. 

Though we had traveled east, and therefore further away from both our parents, moving to Ann Arbor did put us closer to family by the airline metric. I appreciated being able to fly back to Portland, Oregon on a direct flight, rather than connecting through O'Hare and potentially being trapped by another snowstorm or thunderstorm.

\section{The job and its origin story}

%and what the actual job entails

The founder of Math Reviews, Otto Neugebauer, was a German historian of mathematics. He interpreted the mathematical content of Babylonian tablets and worked with Richard Courant and Ferdinand Springer to found the famous ``yellow series'' of mathematical textbooks. In 1931, Courant, Springer, Neugebauer, and Harald Bohr created \emph{Zentralblatt f\"{u}r Mathematik}, a journal that published information about new research mathematics, with summaries or \emph{reviews} of each paper contributed by volunteer mathematicians. As the Nazis came to power, they forced the Jewish mathematician Tullio Levi-Civita off the \emph{Zentralblatt} editorial board and insisted that German mathematics should only be reviewed in German. Neugebauer was intellectually and personally committed to the idea that mathematics prospers through multicultural collaboration. He resigned from the \emph{Zentralblatt} editorial board, came to the United States, and founded a new review journal in coordination with the American Mathematical Society. The first issue of \emph{Mathematical Reviews} was published in January 1940.

Today, Math Reviews is a database (MathSciNet) containing more than four million pieces of mathematical research. As an editor, I'm responsible for looking at research, determining whether it should be part of the database, classifying the mathematics involved, matching research to potential reviewers, and editing reviews so that they make mathematical sense. There are sixteen mathematical editors at MR, and between us we split up all of research mathematics, dividing it by Math Subject Classification number. My current responsibilities, as of the beginning of 2021, include a mix of algebra and geometry, as well as math history:

\begin{description}
\item[01] History of math
\item[12] Field theory
\item[13] Commutative algebra
\item[14] Algebraic geometry
\item[51] Geometry (axiomatic treatments)
\item[52] Convex and discrete geometry
\item[55] Algebraic topology.
\end{description}

When I accepted the position at Math Reviews, the editing part of the job loomed largest in my mind. I spend lots of time thinking about mathematics (especially when a reviewer finds a potential error in a paper!) The job also draws heavily on language skills: because mathematics research is international, effective editing often involves identifying an implicit translation. 

As I settled into my role as an editor, I realized how much I was learning about the publishing industry. I've seen the whole spectrum of research journals, from insightful and excellent all the way down to predatory. I know more about the ins and outs of Open Access schemes and publication models, and have a far better sense of how to identify the right fit for a particular paper. 

The work is important to me for bigger, more symbolic reasons as well. I'm a person who is interested in networks, webs, and complexities. Making the connections between different kinds of mathematics visible is a task I find compelling.

\section{Writing, editing, and self-definition}

% yearly evaluations with actual raises... but much more about balancing workload, missing the structure of tenure/promotion
% established research agenda - that part was humming along
%
%* Returning to creative writing, and what I learned about my own identity in the process

Though my path to tenure was relatively easy, the trinity of teaching, research, and service still dominated all the ways I used my time. As a mathematical editor, it's important for me to stay engaged in the research community, but the core ways my job is evaluated are much more closely related to the day-to-day work of editing. That gave me new freedom to chart my personal goals. 

I was confident in my research program: my physics and number theory collaboration was humming along. But I also wanted to rekindle my writing. I have always written poetry and fiction as well as math---I even talked my way into an MFA poetry course as a graduate student---but as a new teacher, I'd focused my empathetic imagination on my students and their approaches to mathematics. In a more solitary job, I found my imagination stretching in different ways.

As I wrote more, I made more writerly friends, and reawakened some old friendships. I talked my friend Marie Vibbert through the intricacies of theorem naming for a story, ``Loitering with Mathematical Intent,'' about a grad student struggling with research block in a post-capitalist society. Like Marie, my new friends were mostly writing science fiction. They were interested in robots and spaceships, but also in the ways our society constructs norms around gender and sexuality, and how future societies might form different expectations.

I've always written from the points of view of people with all sorts of genders and sexualities. I tended to view this choice as arbitrary, akin to deciding whether to multiply matrices on the right or the left. However, as I wrote more, I recognized that certain kinds of love stories---including stories about two furious, ambitious women---were genuinely easier for me to tell. 

This personal recognition came with a shift in personal labels, to ``bisexual'', or sometimes the more all-encompassing ``queer''. It had been a long time since I'd considered how to classify my sexuality. As a teenager in the '90s, I had the idea that there was a continuum from straight to gay, with bisexuality exactly in the middle, and if you didn't quite fit at one of those points, you rounded. I liked dating boys, so ``rounding down to straight'' made sense to me. My occasional bouts of intense admiration for particular women didn't fit into a narrative that I recognized, so they didn't affect my self-classification scheme.

But the '90s were a long time ago. Many of the models of queerness I had as a teenager seem irrelevant or offensive now. I remember watching \emph{Chasing Amy}, an entire movie dedicated to debating the question of whether bisexuality exists, in the theater! 

Lots of people around my age are reassessing how they identify in light of new information. Personally, I had the intense and mixed emotions that sometimes arise after struggling with a proof, where you're torn between wanting people to share in your revelation and being embarrassed you didn't notice the trick long ago. But there was also a comfort in reclassification: my friend groups have always been very queer, in the more modern, ecumenical sense, and being able to describe that commonality was satisfying.

I picked up another useful strategy from hanging around with lots of gender nerds: asking yourself what parts of gender feel good. I have hard-won mental radar systems for when people are trying to box me in with gendered expectations, particularly professional ones. I know from experience that when I'm in a room full of women in mathematics, some of those expectations lift, which does feel good. But women in STEM networking events have also taught me that many women are stressed almost to breaking by the way society insists on a conflict between math and femininity---I'm thinking, for example, of my undergraduate physics lab partner, who loved pink cardigans and shopping at the Gap, had a perfect lab notebook, did her PhD at Cornell, and now works for a consulting firm. In contrast, I've tended to feel as if I were deflecting a bunch of attacks through sheer obliviousness (and occasionally tripping over javelins that I hadn't known were in the air). 

I noticed, once I started paying more attention, that I actively like having a gender that confuses people. To pick a small example, the Twitter algorithm has inferred, based on the accounts I follow and the tweets I interact with, that I'm a man. This is not correct. I could have changed the setting, but messing up the algorithm made me happy, so I let it be. (When I checked, Twitter also believed that I was interested in education, travel, coffee, and football. The football interest was entirely due to my friend Marie's dedication to the Cleveland Browns; I deleted it.) 

Here's a larger example. I could have responded to my middle school's pressures, all those years ago, by attempting to excel at the skills of femininity. I would never have been popular, but I could have learned how to wear makeup and sorted out which bands I was supposed to admire. Instead, I identified the absolute minimum effort required to avoid harassment, and stopped there. My favorite article of clothing as a seventh-grader was a soft blue hoodie; my favorite article of clothing as a mid-career mathematician is a brown one. Admitting my preferences were personal, real, and consistent made it a lot easier for me to appreciate other people's impulses toward overt femininity without feeling as if I was boxed in by norms, or worrying that I was betraying feminism by hating pink too much.

This is the point, in a conventional twenty-first century American coming out narrative, where I should tell you what my pronouns are. The problem is that I have a distinct and genuine emotional response to stating my pronouns, and that response is fury that I'm being asked to be complicit in a broken system. I was fine using she/her pronouns, as long as I didn't have to pretend that I \emph{liked} them. 

I'm happy with they/them pronouns in contexts where they feel legitimately amorphous. Many kind and well-meaning people are tempted to replace the binary gender system by a triple gender system, and the change from $n$ boxes to $n+1$ boxes does not bring me joy. For that matter, I'm fine with he/him pronouns, if they're not attached to the assumption that all serious mathematicians use he/him. 

In practice, I pick some list of pronouns in some order based on the effect I want to have on the social group I'm joining; my current default for mathematical settings is she/they. I feel a little bit like I'm lying, no matter which list I choose, but it's the kind of lie you tell when you're teaching freshman calculus: sometimes you have to start with intuition and examples, rather than precise definitions.

\section{The once and future Fields Institute}

%
%* Organizing the Homological Algebra and Mirror Symmetry workshop at Fields, and recognizing myself as a role model. 

I met Tyler Kelly at the Fields Institute in 2011. We were both attending a conference associated with a special semester on mirror symmetry. I was getting ready to start my job at UWEC, and Tyler was a graduate student at the University of Pennsylvania. Tyler was a bit of a rising star: they had already done research in mirror symmetry as an undergrad. For Tyler, graduate school was also a way to escape from their home state of Georgia. They were keenly aware of the prejudice they faced as a queer person, both within and outside the mathematical community.

Tyler's undergraduate project involved computing Picard--Fuchs equations, something I had worked on as a grad student. Picard--Fuchs equations are differential equations that describe the way an integral of a holomorphic form changes as you move through a family of geometric spaces. It's a weird, hybrid subject, hovering between intense combinatorial computation and an inspiration that feels very analytic. When I was a graduate student, my algebro-geometric peers were working on high-powered algebraic questions, leveraging subtle properties of rings. I wasn't sure how these squishier equations fit into the field. But Tyler loved them! When I formed a new collaboration to investigate the interface between mirror symmetry and number theory, writing grants to visit the Banff International Research Station and the American Institute of Mathematics, I knew Picard--Fuchs equations would be involved, in their new-to-me arithmetic guise. I invited Tyler.

Eight years and several major papers later, I was one of the co-organizers for a new special semester at Fields, this one on homological mirror symmetry. Meanwhile, Tyler had started a job as a professor at the University of Birmingham, in England. We teamed up to run a summer school for undergraduates on homological mirror symmetry, inviting Amanda Francis, one of my colleagues at Math Reviews, to join us as a speaker. 

Homological mirror symmetry is huge and intimidating; it's rare to approach it at the undergraduate level. But we saw a way in, using a framework called Berglund--H\"{u}bsch--Krawitz mirror symmetry. We recruited an enthusiastic, insightful, accomplished, collaborative group of American and Canadian students. They rose to the challenge joyfully, digging into the mathematics and physics. 

One morning, before our first lecture, a student mentioned that she had read one of my stories and become curious about the very short form of fiction known as flash. I knew which story she meant: it was something I'd published with the Lesbian Historic Motif Podcast, riffing on a medieval Japanese novel where a woman marries another woman. I wasn't used to that kind of contextual collapse. But this was a context where being casually queer and intensely mathematical, simultaneously, mattered. As organizers, we were the grownups; we were demonstrating what mathematical community could be. 

I took the lesson to heart: bridging queer and mathematical communities matters. I was already thinking carefully about how to ask about gender and pronouns in official settings, and when it's important to share queer history in math contexts. Queer identity wove its way into my ongoing professional projects in small but tangible ways. For example, as a Feature Columnist for the AMS, I wrote about the trans math-major heroine of Courtney Milan's novel \emph{Hold Me} and related real-world research; in connection with ongoing efforts to broaden Wikipedia's perspective, I made a Wikipedia page for the LGBTQ+ math association Spectra. I've been talking to people in my professional societies about how mathematicians can oppose homophobia and transphobia. And though it's awkward---though it would be easier if I fit into a tidy narrative---I've been working, in situations where I'm acting as a role model in part because of my gender, to be transparent about how complicated my relationship to my own gender is.

\end{document}